\documentclass[10pt,sort&compress]{elsarticle}
\usepackage{amsfonts,amssymb,amsmath,latexsym,indentfirst}
\usepackage{wasysym,mathrsfs,enumerate,fontenc}
\usepackage{setspace}
\usepackage[left=3cm,top=3cm,right=3cm,bottom=3cm]{geometry}

\numberwithin{equation}{section}

\newtheorem{lemma}{Lemma}[section]
\newtheorem{Theorem}{Theorem}[section]
\newtheorem{corollary}{Corollary}[section]

\newcommand{\R}{\ensuremath{{\mathbb{R}} }}

\newcommand{\fim}{\hfill$\square$\\ \\}
\newcommand{\proof}{\noindent\textbf{Proof.}\quad}

\begin{document}
\begin{frontmatter}
\title{Asymptotic Behavior of Solutions for the Cauchy Problem of a Dissipative Boussinesq-Type Equation}
\journal{Journal}
 \author{\textbf{Amin Esfahani}\corref{cor1}}

 \address{School of Mathematics and Computer Science,  Damghan University,
 	Damghan 36715-364, Iran.\\ E-mail:
 	amin@impa.br, esfahani@du.ac.ir}
 \cortext[cor1]{Corresponding author}
 \author{\textbf{Hamideh B. Mohammadi}}
 \address{School of Mathematics and Computer Science,  Damghan University,
 	Damghan 36715-364, Iran.}

\begin{abstract}
We consider the Cauchy problem for an evolution equation modeling bidirectional surface waves in a convecting fluid. 
Under small condition on the initial value, the existence and
asymptotic behavior of global solutions in some time
weighted   spaces
are established by the contraction mapping principle.
\end{abstract}
\begin{keyword}
Dissipative Boussinesq equation \sep Well-posedness\sep Sobolev spaces.\\
\emph{MSC[2010]}: 35B30\sep 35Q55\sep 35Q72.
\end{keyword}
\end{frontmatter}

\section{Introduction}
The Kuramoto-Sivashinsky (KS) equation
\begin{equation}
u_t+\gamma u_{xxxx}+\alpha u_{xx}+uu_x=0
\end{equation}
is a well-known model of one-dimensional turbulence derived in various physical contexts such as chemical-reaction waves, propagation
of combustion fronts in gases, surface waves in a film of a viscous liquid flowing along an inclined
plane, patterns in thermal convection, rapid solidification (see e.g. \cite{hm,sivash,yk}), where $\alpha$ and $\gamma$ are constant coefficients accounting for the long-wave instability (gain) and short-wave dissipation, respectively. By combining the dispersive effects of the KdV equation and the dissipative effects of the KS equation, the Kuramoto-Sivashinsky-Korteweg-de Vries (KS-KdV) equation
\begin{equation}
u_t+u_{xxx}+\gamma u_{xxxx}+\alpha u_{xx}+uu_x=0
\end{equation}
appears; which  was first introduced by Benney \cite{benney}. This equation
finds various applications   in the study of unstable drift
waves in plasmas \cite{cktr}, fluid flow along an inclined plane \cite{benney,tk} convection in fluids with a free surface \cite{bd,bd-1,ad,gv} the Eckhaus
instability of traveling waves \cite{jpbcrk}, in solar
dynamo wave \cite{mk}, hydrodynamics and other fields \cite{cv,eig,oe}.

The derivation of this equation in the physical situations
mentioned above involves the assumption of unidirectional
waves.  The assumption of unidirectional waves for surface waves was removed in  \cite{kkpm,nv} and a modified Boussinesq system of equations was derived. One of these type of equation is the following dissipative Boussinesq equation:
\begin{equation} \label{equation 11}
u_{tt} - \Delta u + \Delta^2 u + \alpha\Delta u_t +
\gamma\Delta ^2u_t = \Delta (\beta f(u_t) + g(u)).
\end{equation}  
Here $u= u(x,t)$ is the unknown function of $x= (x_1,\cdots , x_n)
\in \R^n$ $, t>0 $ and $\beta >0$ and $\alpha\in\mathbb{R}$ are
constants. The term  $u_t$ represents a frictional function
dissipation, and the nonlinear term $f(v)$ and $g (v)$ are smooth
functions of $v$ under considerations and satisfies $f(v)= O(\vert
v \vert ^2)$ and $g(v)= O(\vert v \vert ^2)$ for $v\rightarrow 0$. Equation \eqref{equation 11} arises in the  study of the stability of one-dimensional periodic patterns in
 systems with Galilean invariance and also the oscillations of elastic
 beams \cite{christov,cmp,cv-2}.  Ignoring the dissipation, \eqref{equation 11} turns into the classical Boussinesq equation
\begin{equation} \label{bouss}
u_{tt} -\Delta u \pm
\Delta ^2u = \Delta (u^2);
\end{equation}
appeared not only in the study of the dynamics of thin inviscid layers with free surface but also in the
study of the nonlinear string, the shape-memory alloys, the propagation of waves in elastic rods and
in the continuum limit of lattice dynamics or coupled electrical circuit. When $\gamma=\beta=0$, the existence, uniqueness and long-time asymptotic of solutions to the Cauchy problem and the initial
boundary value problem of equation \eqref{equation 11}  has been studied by several authors, see for instance \cite{esfahani,liu,var-1,var-2,var-3,var-4,var-5,var-6,wang-0,wang} and references therein.

In this paper we study the asymptotic behavior of solutions of the Cauchy problem associated to \eqref{equation 11} 
with the initial values
\begin{equation} \label{equation 12}
u(0) = u_0(x),\ \ \ u_t(0)= u_1(x).
\end{equation}

The article is organized as follows. In Section \ref{section2} we obtain the solution formula
of \eqref{equation 11} and study the decay property of the solution operators
appearing in the solution formula. Then, in Section \ref{section3}, we discuss the
linear problem and show the decay estimates of the solutions in $L^1$. We prove global existence and asymptotic behavior of solutions for the Cauchy problem \eqref{equation 11} and \eqref{equation 12} in $L^2$ in Section
\ref{section4}.

Throughout this paper we assume   $\gamma=1\leq-\alpha$.
\section{Decay property of the linear part}\label{section2}
The aim of this section is to derive the solution formula for the
problem \eqref{equation 11} and \eqref{equation 12}. First of all,
we investigate the linear equation of \eqref{equation 11}.
\begin{equation} \label{equation 21}
	u_{tt} - \Delta u + \Delta^2 u + \alpha\Delta u_t +
	\Delta ^2u_t =0
\end{equation}
with the initial data \eqref{equation 12}.
\\By applying the Fourier transform to \eqref{equation 21}  we have
\begin{equation} \label{equation 22}
	\hat{u}_{tt} + (\vert\xi\vert^4 - \alpha \vert \xi\vert
	^2)\hat{u}_t+ (\vert\xi\vert^2 +\vert \xi\vert^4)\hat{u}
	= 0.
\end{equation}
The corresponding initial values are given as
\begin{equation} \label{equation 23}
	t=0:\ \ \hat{u}=\hat{u}_0(\xi),\ \ \ \ \hat{u}_t=
	\hat{u}_1(\xi).
\end{equation}
The characteristic equation of \eqref{equation 22} is
\begin{equation}
	\label{equation 24}
	\lambda^2+ (\vert\xi\vert ^4- \alpha\vert\xi\vert^2)\lambda+
	(\vert\xi\vert^4+|\xi|^2)=0.
\end{equation}
Let $\lambda=\lambda_{\pm} (\xi)$ be the corresponding
eigenvalues, i.e
\begin{equation} \label{equation 25}
	\lambda_{\pm}(\xi)= \frac {(\alpha\vert\xi\vert^2-
		\vert\xi\vert^4)\pm
		\sqrt{(\vert\xi\vert^4-\alpha\vert\xi\vert^2)^2-4(\vert\xi\vert^2+\vert\xi\vert^4)}}{2}
	.\end{equation}
The solution to the problem \eqref{equation 22} and \eqref{equation 23} is given in the form
\begin{equation}
	\label{equation 26}
	\hat{u}(\xi,t)= \hat G(\xi,t) \hat{u}_1(\xi)+\hat{H}(\xi,t)
	\hat{u}_0(\xi),
\end{equation}
where
\begin{equation}
	\label{equation 27}
	\hat{G}(\xi,t)= \frac{1}{\lambda_+(\xi) - \lambda_-(\xi)}(e^{\lambda_+(\xi)t}- e^{\lambda_-(\xi)t})
\end{equation}
and
\begin{equation}
	\label{equation 28}
	\hat{H}(\xi,t)= \frac{1}{\lambda_+(\xi) - \lambda_-(\xi)}(\lambda_+(\xi)e^{\lambda_-(\xi)t}-
	\lambda_-(\xi)e^{\lambda_+(\xi)t}).
\end{equation}
Let
\begin{equation}
	\label{equation 29}
	G(x,t)= F^{-1}[\hat{G}(\xi,t)](x)
\end{equation}
and
\begin{equation} \label{equation 30}
	H(x,t)= F^{-1}[\hat{H}(\xi,t)](x),
\end{equation}
where $F^{-1}$ denotes the inverse Fourier transform. With applying
$F^{-1}$ to \eqref{equation 26}, we obtain
\begin{equation}
	\label{equation 10}
	u(t) =G(t)\ast u_1+ H(t)\ast u_0.
\end{equation}
By the Duhamel principle, we obtain the solution formula to
\eqref{equation 11} and \eqref{equation 12}
\begin{equation}
	\label{equation 1}
	u(t)=G(t)\ast u_1+ H(t)\ast u_0+ \int_0^t G(t-\tau)\ast \Delta(f(u (\tau))+\beta
	g(u_t))(\tau)d\tau.
\end{equation}

Now we study the decay property of the linear
equation \eqref{equation 11}. Our aim is to prove the following
decay estimates of the solution operators $G(t)$ and $H(t)$
appearing in (\ref{equation 10})

\begin{lemma}
	The solution of \eqref{equation 22} and \eqref{equation 23}
	satisfies
	\begin{equation}\label{equation 31}
		\vert\xi\vert^2(1+\vert\xi\vert^2)\vert\hat{u}(\xi,t)\vert^2+
		\vert\hat{u}_t(\xi,t)\vert^2\leq Ce^{- c\omega(\xi)t}(\vert\xi\vert^2(1+\vert\xi\vert^2)\vert\hat{u}_0(\xi)\vert^2+ \vert\hat{u}_1(\xi)\vert^2)
	\end{equation}
	for $ \xi\in R^n$ and $t\geq 0$, where $\omega(\xi)= \vert\xi\vert^2/(1+|\xi|^2).$
\end{lemma}
\proof
	By multiplying \eqref{equation 22} by $\bar{\hat{u_t}}$
	and taking
	the real part ,  we deduce that
	\begin{equation}\label{equation32}
		\frac{1}{2}\frac{d}{dt}(|\hat
		{u}_t|^2+(|\xi|^2+|\xi|^4)|\hat{u}|^2)+(|\xi|^4-\alpha|\xi|^2)|\hat{u}_t|^2=0
	\end{equation}
	Multiplying \eqref{equation 22} by $\bar{\hat{u}}$ and take the real part
	yields
	\begin{equation}\label{equation33}
		\frac{1}{2}\frac{d}{dt}((|\xi|^4-\alpha|\xi|^2)|\hat{u}|^2+ 2 Re
		(\hat{u}_t.\bar{\hat{u}}))+ (|\xi|^2+
		|\xi|^4|)|\hat{u}|^2-|\hat{u}_t|^2=0
	\end{equation}
	Multiplying both sides of \eqref{equation32} and \eqref{equation33} by $(1+|\xi|^2)$ and
	$|\xi|^2$ respectively, summing up the products yields
	\begin{equation}\label{equation 34}
		\frac{d}{dt} E+F=0,
	\end{equation}
	where
	\begin{equation}
		\nonumber  E=(1+|\xi|^2)|\hat{u_t}|^2+\{(1+|\xi|^2)(|\xi|^2+
		|\xi|^4)+ |\xi|^2 (|\xi|^4- \alpha|\xi|^2)\}|\hat{u}|^2+ 2|\xi|^2
		Re(\hat{u_t}\bar{\hat{u}})
	\end{equation}
	and
	\begin{equation}
		\nonumber  F=\{2(1+|\xi|^2)(|\xi|^4-\alpha|\xi|^2)-2|\xi|^2\}\  |\hat{u_t}|^2+
		2|\xi|^2 (|\xi|^2+|\xi|^4)|\hat u|^2.
	\end{equation}
	It is easy to see that
	
	\begin{equation}\label{equation 35}
		C(1+|\xi|^2)E_0\leq E\leq C(1+|\xi|^2)E_0,
	\end{equation}
	where
	\begin{equation}
		\nonumber  E_0= |\hat{u_t}|^2+ |\xi|^2(1+|\xi|^2)|u|^2.
	\end{equation}
	Noting that $F\geq|\xi|^2 E_{0}$ and with \eqref{equation 35}, we
	obtain
	\begin{equation}\label{equation 36}
		F\geq c \ \omega(\xi)E,
	\end{equation}
	where
	\begin{equation}
		\nonumber \omega(\xi)= \frac{|\xi|^2}{1+|\xi|^2}.
	\end{equation}
	Using \eqref{equation 34}and \eqref{equation 36}, we get
	\begin{equation}
		\nonumber  \frac{d}{dt}E +c
		\ \omega(\xi)E\leq 0.
	\end{equation}
	Thus
	\begin{equation}
		\nonumber E(\xi,t)\leq e^{-c \ \omega(\xi)t}E(\xi,0),
	\end{equation}
	which together with \eqref{equation 35} proves the desired
	estimate \eqref{equation 31}.
\fim
\begin{lemma}
	Assume that $\hat{G}(\xi,t)$ and $\hat{H}(\xi,t)$ are fundamental
	solutions of \eqref{equation 21} in the Fourier space, which are
	given explicitly in \eqref{equation 27} and \eqref{equation
		28}.Then we have the pointwise estimates
	\begin{equation}
		\label{equation 37}
		|\xi|^2(1 +|\xi|^2)|\hat{G}(\xi,t)|^2+ |\hat{G_t}(\xi,t)|^2\leq
		C e^{-c\omega(\xi)t}
	\end{equation}
	and
	\begin{equation}
		\label{equation 38}
		|\xi|^2(1 +|\xi|^2)|\hat{H}(\xi,t)|^2+
		|\hat{H_t}(\xi,t)|^2\leq
		C |\xi|^2 (1+|\xi|^2)e^{-c\omega(\xi)t},
	\end{equation}
	
	for $\xi \in \mathbb{R}^n$ and $t\geq0$, where
	$\omega(\xi)=\frac{|\xi|^2}{1+|\xi|^2}$.
\end{lemma}
\proof
	If $\hat{u}_0(\xi)=0$ , then from \eqref{equation 26} we get
	\begin{equation}
		\nonumber \hat{u}(\xi,t)= \hat{G}(\xi,t) \hat{u_1}(\xi),\ \ \
		\hat{u}_t(\xi,t)= \hat{G}_t(\xi,t)\hat{u}_1(\xi).
	\end{equation}
	Substituting the equalities into \eqref{equation 31} with
	$\hat{u}_0(\xi)=0$ we
	obtain \eqref{equation 37}.
	In what follows, we consider $\hat{u}_1(\xi)=0$. We have  from \eqref{equation 26}
	that
	\begin{equation}
		\nonumber \hat{u}(\xi,t)= \hat{H}(\xi,t) \hat{u}_0(\xi),\ \ \
		\hat{u}_t(\xi,t)= \hat{H_t}(\xi,t)\hat{u}_0(\xi).
	\end{equation}
	Substituting the equalities into \eqref{equation 31} with
	$\hat{u}_1(\xi)=0$, we obtain \eqref{equation 38}  , which
	together with \eqref{equation 37}, we have completed the proof of
	the lemma.
\fim
\begin{lemma}
	Let $l, k , j$  be nonnegative integers and assume that $1\leq p\leq
	2$. Then we have
	\begin{eqnarray}
		\nonumber
		\|\partial_x^k G(t)\ast\phi\|_{L^2} &\leq & C(1+t)^{-\frac{n}{2}(\frac{1}{p}-\frac{1}{2})-
			\frac{k-j}{2}}\ \|\partial_x^j \phi\|_{\dot{W}^{-1 ,p}}\\&&+C
		e^{-ct} \|\partial_x^{k+l-2}\phi\|_{L^2}\label{equation 39},
	\end{eqnarray}
	\begin{eqnarray}
		\nonumber
		\|\partial_x^k H(t)\ast\psi\|_{L^2}& \leq & C(1+t)^{-\frac{n}{2}(\frac{1}{p}-\frac{1}{2})-
			\frac{k-j}{2}}\ \|\partial_x^j \psi\|_{L^p}\\&& +C
		e^{-ct} \|\partial_x^{k+l}\phi\|_{L^2}\label{equation 310},
	\end{eqnarray}
	for $0\leq j\leq k$, where $k+l-2\geq0$ in \eqref{equation 39}. Similarly, we have
	\begin{eqnarray}
		\nonumber
		\|\partial_x^k G_t(t)\ast\phi\|_{L^2}& \leq & C(1+t)^{-\frac{n}{2}(\frac{1}{p}-\frac{1}{2})-
			\frac{k+1-j}{2}}\ \|\partial_x^j \phi\|_{\dot{W}^{-1 ,p}}\\& &+C
		e^{-ct} \|\partial_x^{k+l}\phi\|_{L^2}\label{equation 311},
	\end{eqnarray}
	\begin{eqnarray}
		\nonumber
		\|\partial_x^k H_t(t)\ast\psi\|_{L^2}& \leq& C(1+t)^{-\frac{n}{2}(\frac{1}{p}-\frac{1}{2})-
			\frac{k+1-j}{2}}\ \|\partial_x^j \psi\|_{L^p}\\ &&+C
		e^{-ct} \|\partial_x^{k+l+2}\phi\|_{L^2}\label{equation 312},
	\end{eqnarray}
	for $0\leq j\leq k+1$.

\end{lemma}
\proof
	We only give a proof of \eqref{equation 39}. We apply the
	Plancherel theorem and use the pointwise estimate for $\hat{G}$ in
	\eqref{equation 37}. This gives
	\begin{eqnarray}
		\nonumber \|\partial_x^k G _t(t)\ast\phi\|_{L^2}^2 &=&\int_{\mathbb{R}^n}|\xi|^{2k}|
		\hat{G}(\xi,t)|^2\|\hat{\phi}(\xi)|^2
		d\xi\\
		&=&  \int_{|\xi|\leq
			1}|\xi|^{2k}|\hat{G}(\xi,t)|^2\hat{\phi}(\xi)|^2d \xi\nonumber
		+\int_{|\xi|\geq
			1}|\xi|^{2k}|\hat{G}(\xi,t)|^2\hat{\phi}(\xi)|^2 d\xi\nonumber\\
		& \nonumber \leq& \int_{|\xi| \leq 1}|\xi|^{2k-2}e^{- c |\xi|^2
			t}|\hat{\phi}(\xi)|^2 d\xi  \\
		&\nonumber &+ C\int_{|\xi|\geq 1}e^{-c \omega (\xi)t}
		|\xi|^{2k}(|\xi|^2(1+|\xi|^2))^{-1}|\hat{\phi}(\xi)|^2
		d\xi\\
		& \nonumber \leq& C ||\ |\xi|^{j-1}\hat{\phi}(\xi)\
		\|_{L^{\acute{p}}}^2 \ (\int_{|\xi|\leq 1}|\xi|^{2(k-j)q} e^{-cq|\xi|^2
			t}
		d\xi)^{\frac{1}{q}}
		\\
		& \nonumber &+ C e^{-ct} \int_{|\xi| \geq 1}|\xi|^{2k-4}
		|\hat{\phi}(\xi)|^2d\xi\\
		& \nonumber \leq& C ||\ |\xi|^{j-1}\hat{\phi}(\xi)\
		\|_{L^{\acute{p}}}^2 \ (||\ |\xi|^{2(k-j)} e^{-c|\xi|^2
			t}||_{L^q}\\
		& \nonumber &+ C e^{-ct} \int_{|\xi| \geq
			1}|\xi|^{2(k+l-2)}|\hat{\phi}(\xi)|^2 d\xi,
	\end{eqnarray}
		where we used H\"{o}lder inequality with $\frac{1}{q}+
	\frac{2}{\acute{p}}=1$, $\frac{1}{p}+\frac{1}{\acute{p}}=1$.
	With a straight computation, we obtain
	\begin{equation}
		\nonumber \||\xi|^{2(k-j)} e^{-c|\xi|^2
			t}\|_{L^q(|\xi|\leq1)}\leq C
		(1+t)^{-n(\frac{1}{p}- \frac{1}{2})- (k-j)}.
	\end{equation}
	It follows from the Hausdorff-Young inequality that
	\begin{equation}
			\nonumber \||\xi|^{j-1}\hat{\phi}(\xi)\|_{L^{\acute{p}}}\leq
		\|\partial_x^j\phi\|_{\dot{W}^{-1,p}}
	\end{equation}
	Combining the above three inequalities  yields \eqref{equation
		39}. Similarly, we can prove \eqref{equation 310}-\eqref{equation 312}. Thus the lemma is proved.
\fim
Immediately we have from previous lemma the following corollary.

\begin{corollary}
	Let $1\leq p\leq2$,and let k ,j and l be nonnegative integers.
	Also, assume that $G(x,t)$ and $H(x,t)$ be the fundamental
	solution of \eqref{equation 21} which are given in \eqref{equation 27} and \eqref{equation 28}, respectively. Then we have
	\begin{eqnarray}
		\nonumber
		\|\partial_x^k G(t)\ast\Delta g\|_{L^2}& \leq & C(1+t)^{-\frac{n}{2}(\frac{1}{p}-\frac{1}{2})-
			\frac{k+1-j}{2}}\ \|\partial_x^j g\|_{L^p} \\&&+C
		e^{-ct} \|\partial_x^{k+l}g\|_{L^2}\label{equation 313},
	\end{eqnarray}
	for $0\leq k\leq j+1$. It also for $0\leq k\leq j+2$ holds that
	\begin{equation}\label{equation 314}
			\|\partial_x^k G_t(t)\ast\Delta g\|_{L^2} \leq C(1+t)^{-\frac{n}{2}(\frac{1}{p}-\frac{1}{2})-
			\frac{k+2-j}{2}}\ \|\partial_x^j g\|_{L^p} +C
		e^{-ct} \|\partial_x^{k+l+2}g\|_{L^2}.
	\end{equation}
\end{corollary}

\section{Global existence and asymptotic behavior of solutions for $L^1$}\label{section3}

The aim of this section is to prove the existence and asymptotic
behavior of solutions to \eqref{equation 11} and \eqref{equation 12} with $L^1$ data. We first state the
following lemma, which comes from \cite{zheng}.
\begin{lemma}
	Assume that $f=f(v)$ is smooth function, where
	$v=(v_1,\ldots,v_n)$ is a vector function. Suppose that $f(v)=
	O(|v|^{1+\theta})(\theta\geq 1$ is an integer) when $|v|\leq v_0$.
	Then, for the integer $m\geq 0$, if  $v,w \in
	W^{m,q}(\mathbb{R}^n)\bigcap L^p(\mathbb{R}^n)\bigcap
	L^\infty(\mathbb{R}^n)$ and $\|v\|_{L^{\infty}}\leq
	v_0$,$\|w\|_{L^{\infty}}\leq v_0$, then $f(v)-f(w) \in
	W^{m,r}(\mathbb{R}^n)$. Furthermore, the following inequalities
	hold:
	\begin{equation}
		\label{equation 41}
		\|\partial_x^m f(v)\|_{L^r} \leq  C\|v\|_{L^p}\|\partial_x^m
		v\|_{L^q}\|v\|_{L^\infty}^{\theta-1}
	\end{equation}
	and
	\begin{eqnarray}
		\label{equation 42}
		\|\partial_x^m( f(v)-f(w))\|_{L^r} & \leq & C\{(\|\partial_x^m v\|_{L^q})
		\|v-w\|_{L^p}+\\&&(\|v\|_{L^p}+
		\|w\|_{L^p}\|\partial_x^m(v-w)\|_{L^q})\}(\|v\|_{L^{\infty}}+\|w\|_{L^\infty})^{\theta-1}\nonumber,
	\end{eqnarray}
	where $\frac{1}{r}=\frac{1}{p}+\frac{1}{q}$ and $1\leq p,q,r\leq
	+\infty$.
\end{lemma}
Based on the decay estimates of solutions to the linear problem
\eqref{equation 21},we define the following solution space:
\begin{equation}
	\nonumber X=\{u \in C([0,\infty);H^{s+2}(\mathbb{R}^n))\bigcap
	C^1([0,\infty);H^s(\mathbb{R}^n)):\|u\|_X<\infty\},
\end{equation}
where
\begin{equation}
	\nonumber \|u\|_X= \sup_{t\geq 0} \left \{\sum_{k\leq
		s+2}(1+t)^{\frac{n}{4}+\frac{k}{2}}\|\partial_x^k u(t)\|_{L^2}+
	\sum_{k\leq s}(1+t)^{\frac{n}{4}+\frac{k}{2}}\|\partial_x^k
	u_t(t)\|_{L^2}\right\}.
\end{equation}
For $ R>0 $, we define
$$X_R= \{u\in X:\|u\|_X\leq R\}.$$
\begin{Theorem}\label{thm 41}
Let $n\geq 1,s\geq
\max\{0,[\frac{n}{2}]-1\}$.	Suppose that $u_0 \in H^{s+2}(\mathbb{R}^n)\bigcap
	L^1(\mathbb{R}^n)$, $u_1 \in H^s(\mathbb{R}^n)\bigcap
	\dot{W}^{-1,1}(\mathbb{R}^n)$ and $f(v),g(v)$ are smooth and
	satisfies $f(v)= O(v^2),\ g(v)=O(v^2)$ for $v\rightarrow 0$. Put
	\begin{equation}
		\nonumber E_0:= \|u_0\|_{L^1}+
		\|u_1\|_{\dot{W}^{-1,1}}+\|u_0\|_{H^{s+2}}+ \|u_1\|_{H^s}.
	\end{equation}
		If $E_0$ is suitably small,the Cauchy problem \eqref{equation 11}
	and \eqref{equation 12}has a unique global solution $u(x,t)$
	satisfying
	\begin{equation}
		\nonumber X=u \in C([0,\infty);H^{s+2}(\mathbb{R}^n))\bigcap
		C^1([0,\infty);H^s(\mathbb{R}^n)).
	\end{equation}
	Also, the solution satisfies the decay estimate
	\begin{equation}
		\label{equation 43}
		\|\partial_x^k u(t)\|_{L^2}\leq C E_0(1+t)^{-\frac{n}{4}-\frac{k}{2}}
	\end{equation}
	and
	\begin{equation}
		\label{equation 44}
		\|\partial_x^l u_t(t)\|_{L^2}\leq C E_0(1+t)^{-\frac{n}{4}-\frac{l+1}{2}}
	\end{equation}
	for $0\leq k\leq s+2$ and $0\leq l \leq s.$
	
\end{Theorem}
\proof
	The Gagliardo-Nirenberg inequality gives
	\begin{equation}
		\label{equation 45}
		\|u(t)\|_{L^\infty}\leq C \|\partial _x^{s_0}u\|_{L^2}^{\theta}
		\|u\|_{L^2}^{1-\theta}\leq C (1+t)^{- \frac{n}{2}}\|u\|_X
	\end{equation}
	where $s_0= {\frac{n}{2}}+1$, $\theta= \frac{n}{2s_0}$; i.e, $s\geq
	[\frac{n}{2}]-1$. We define
	\begin{equation}
		\nonumber
		\Phi(u)=G(t)\ast u_1+ H(t)\ast u_0 +\int_0^tG(t-\tau)\ast
		\Delta(f(u)-\beta g(u_t))(\tau)d\tau.
	\end{equation}
	We apply $\partial_x^k$ to $\Phi$ and take the $L^2$ norm. We
	obtain
	\begin{equation}\label{equation 46}\begin{split}
			\|\partial_x^k\Phi(u)\|_{L^2}&\leq \|\partial_x^k G(t)\ast
		u_1\|_{L^2}+ \|\partial_x^kH(t)\ast u_0\|_{L^2}\\
		&
	\qquad+C\int_0^t\|\partial_x^kG(t-\tau)\ast
		\Delta(f(u)-\beta g(u_t))(\tau)\|_{L^2}d\tau\\
	& :=I_1+I_2+J
\end{split}	\end{equation}
	First, we estimate $I_1$. We apply \eqref{equation 39} with $p=1$, $j=0$, $l=0$ and get
		\begin{eqnarray}
		\label{equation 47}
		I_1&\leq &C (1+t)^{-\frac{n}{4}-\frac{k}{2}}
		\|u_1\|_{\dot{W}^{-1,1}}+ C
		e^{-ct}\|\partial_x^{(k-2)_+}u_1\|_{L^2}\\
		\nonumber & \leq & C E_0 (1+t)^{- \frac{n}{4}- \frac{k}{2}},
	\end{eqnarray}
	where $(k-2)_+= \max \{k-2,0\}.$

	For the term $I_2$, we apply \eqref{equation 310} with $p=1,j=0$
	and $l=0$. This yields
	\begin{eqnarray}
		\label{equation{48}}
		I_2\leq C (1+t)^{-\frac{n}{4}-\frac{k}{2}}
		\|u_0\|_{L^1}+ C e^{-ct}\|\partial_x^{k}u_0\|_{L^2} \leq C E_0
		(1+t)^{- \frac{n}{4}- \frac{k}{2}}.
	\end{eqnarray}
	Next, we estimate J. Let
	\begin{eqnarray}
		\nonumber
		J&=&\int_0^tG(t-\tau)\ast
		\Delta(f(u)-\beta g(u_t))(\tau)d\tau \\
		\nonumber&=&\int_0^{t/2}G(t-\tau)\ast
		\Delta(f(u)-\beta g(u_t))(\tau)d\tau \\
		\nonumber &&+\int_{t/2}^t G(t-\tau)\ast
		\Delta(f(u)-\beta g(u_t))(\tau)d\tau\\
		\nonumber&=:& J_1+J_2
	\end{eqnarray}
	For the term $J_1$, using \eqref{equation 313} with $p=1,j=0$ and
	$l=0$, we have
	\begin{equation}\label{equation 49}\begin{split}
		J_1&\leq C
		\int_0^{t/2}(1+t-\tau)^{-\frac{n}{4}-\frac{k+1}{2}}\|f(u)(\tau)-\beta
		g(u_t)(\tau)\|_{L^1} d\tau\\
&\qquad+ C
		\int_0^{t/2}e^{-c(t-\tau)}\|\partial_x^k(f(u)-\beta g(u_t))(\tau)\|_{L^2}d\tau\\
	&=:J_{11}+ J_{12}
\end{split}	\end{equation}
	Note that by lemma \eqref{equation 41} we have
	\begin{eqnarray}
		\nonumber \|f(u)\|_{L^1}&\leq & C\|u\|_{L^2}^2 \leq
		CR^2(1+\tau)^{-\frac{n}{2}}\\
		\nonumber|g(u_t)\|_{L^1}&\leq &
		C\|u_t\|_{L^2}^2 \leq CR^2(1+\tau)^{-\frac{n}{2}}.
	\end{eqnarray}
		Therefore we have
	\begin{eqnarray}
		\nonumber
		J_{11}&\leq& CR^2
		\int_0^{t/2}(1+t-\tau)^{-\frac{n}{4}-\frac{k+1}{2}}(1+\tau)^{-\frac{n}{2}}d\tau\\
		\nonumber&\leq&
		CR^2(1+t)^{-\frac{n}{4}-\frac{k+1}{2}}\int_0^{t/2}(1+\tau)^{\frac{n}{2}}d\tau\\
		&\leq&CR^2(1+t)^{-\frac{n}{4}-\frac{k}{2}}\eta(t),\nonumber
	\end{eqnarray}
	where
	\begin{equation}\label{equation 410}
		\eta(t)= \left \{ \begin{array}{c}
			1, \ \ \ \ \ \ \ \ \ \ \ \ \ \ \ \  \ \ \ \ \ \ \ \ \ \ \ \ \   n=1 \\
			(1+t)^{-\frac{1}{2}}\ln(2+t),\ \ \ \ \ \ \ \  n=2\\
			(1+t)^{-\frac{1}{2}},\ \ \ \ \ \  \ \ \ \ \ \ \ \ \ \ \ \  \ \
			n\geq 3.
		\end{array} \right.
	\end{equation}
	We use \eqref{equation 41} and obtain
	
	\begin{equation}
		\label{equation 411}\|\partial_x^k(f(u)(\tau)-\beta
		g(u_t(\tau)))\|_{L^2} \leq
		CR^2(1+t)^{-\frac{n}{4}-\frac{k}{2}-\frac{n}{2}}
	\end{equation}
	Consequently,
	we get
	$$J_{12} \leq CR^2\int_0^{t/2}e^{-c(t-\tau)}(1+\tau)^{-\frac{n}{4}-\frac{k}{2}-\frac{n}{2}}d\tau \leq CR^2e^{-ct}. $$
	Finally, we estimate the term $J_2$ on the time interval
	$[t/2,t]$. Applying \eqref{equation 313} with $ p=2,j=k ,l=0$ and
	using \eqref{equation 411}, we can estimate term $J_2$ as
	\begin{equation}\label{equation 412}\begin{split}
 J_2&\leq C\int_{t/2}^t
		(1+t-\tau)^{-\frac{1}{2}}\|\partial_x^k(f(u)-\beta
		g(u_t)(\tau))\|_{L^2}  d\tau\\
		&\qquad +C\int_{t/2}^t e^{-c(t-\tau)}\|\partial_x^k(f(u)-\beta
		g(u_t)(\tau))\|_{L^2}\\
		&\leq
		CR^2(1+t)^{-\frac{n}{4}-\frac{k}{2}-\frac{n-1}{2}}.
	\end{split}\end{equation}
	Thus we have shown that
	$$J\leq CR^2(1+t)^{-\frac{n}{4}-\frac{k}{2}}\eta(t).$$
	Substituting all these estimates into \eqref{equation 46}, we have
	\begin{equation}
		\label{equation 413} (1+t)^{\frac{n}{4}+\frac{k}{2}}
		\|\partial_x^k\Phi(u)\|\leq CE_0+CR^2,
	\end{equation}
	for $0\leq k \leq s+2.$  It follows from  that \eqref{equation 46}
	\begin{eqnarray}
		\nonumber \Phi(u)_t&=& G_t(t)\ast u_1+ H_t(t)\ast
		u_0\\
		&&+\int_0^tG_t(t-\tau)\ast \Delta(f(u)-\beta
		g(u_t))(\tau)\|_{L^2}d\tau.\label{equation 414}
	\end{eqnarray}
	We use $\partial _x^k$ to $\Phi(u)_t$ and take $L^2$norm.This
	yields
	\begin{equation}\label{equation 415}\begin{split}
		\|\partial_x^k\Phi(u)_t\|_{L^2}&\leq \|\partial_x^k G_t(t)\ast
		u_1\|_{L^2}+ \|\partial_x^kH_t(t)\ast u_0\|_{L^2}\\
		&\qquad
	+C\int_0^t\|\partial_x^kG_t(t-\tau)\ast
		\Delta(f(u)-\beta g(u_t))(\tau)\|_{L^2}d\tau\\
		& =:\acute{I_1}+ \acute{I_2}+\acute{J},
\end{split}	\end{equation}
	for $0\leq k\leq s.$ For the term $\acute{I_1}$, we apply
	\eqref{equation 311} with $p=1,j=0$ and $l=0$ and obtain
	\begin{eqnarray}
		\nonumber
		\acute{I_1}\leq C (1+t)^{-\frac{n}{4}-\frac{k+1}{2}}
		\|u_1\|_{\dot{W}^{-1,1}}+ C e^{-ct}\|\partial_x^{k}u_1\|_{L^2}
		\leq C E_0 (1+t)^{- \frac{n}{4}- \frac{k+1}{2}}.
	\end{eqnarray}
	Also, for the term $\acute{I_2}$, we apply \eqref{equation 312}
	with $p=1,j=0$ and $l=0$ and get
	\begin{eqnarray}
		\nonumber
		\acute{I_2}\leq C (1+t)^{-\frac{n}{4}-\frac{k+1}{2}}
		\|u_0\|_{L^1}+ C e^{-ct}\|\partial_x^{k+2}u_0\|_{L^2} \leq C E_0
		(1+t)^{- \frac{n}{4}- \frac{k+1}{2}}.
	\end{eqnarray}
	To estimate the nonlinear term $\acute{J}$, we divide as
	$\acute{J}=\acute{J_1}+\acute{J_2}$, where $\acute{J_1}$and
	$\acute{J_2}$ correspond to the time intervals $[0,t/2]$ and
	$[t/2,t]$, respectively. By applying \eqref{equation 314} with
	$p=1,j=0$ and $l=0$, we have
	\[\begin{split}
		\acute{J_1}&\leq C
		\int_0^{t/2}(1+t-\tau)^{-\frac{n}{4}-\frac{k+2}{2}}\|f(u)(\tau)-\beta
		g(u_t)(\tau)\|_{L^1} d\tau\\
&\qquad+ C
		\int_0^{t/2}e^{-c(t-\tau)}\|\partial_x^{k+2}(f(u)-\beta g(u_t))(\tau)\|_{L^2}d\tau\\
	&=:\acute{J_{11}}+ \acute{J_{12}}.
	\end{split}\]
	By \eqref{equation 41}, we obtain
	$$\|(f(u)-\beta g(u_t))(\tau)\|_{L^1}\leq
	CR^2(1+\tau)^{-\frac{n}{2}}.$$
	Therefor we get
\[	\begin{split}
			\acute{J_{11}}&\leq CR^2
		\int_0^{t/2}(1+t-\tau)^{-\frac{n}{4}-\frac{k+2}{2}}(1+\tau)^{-\frac{n}{2}}d\tau\\
		&\leq CR^2(1+t)^{-\frac{n}{4}-\frac{k+1}{2}}\eta(t).
	\end{split}\]
	Similarly as before,we can estimate $\acute{J_{12}}$ and obtain
	$\acute{J_{12}}\leq CR^2e^{-ct}.$ Finally, we estimate the term
	$\acute {J_2}$ by using  \eqref{equation 314} with $p=2$, $j=k+2$, $l=0$ and get
	\begin{eqnarray}
		\nonumber\acute{ J_2}&\leq&
		C\int_{t/2}^t \|\partial_x^{k+2}(f(u)-\beta
		g(u_t)(\tau))\|_{L^2}  d\tau\\
		&&\nonumber +C\int_{t/2}^t{-\frac{n}{4}-\frac{k+1}{2}-\frac{n-1}{2}}
		e^{-c(t-\tau)}\|\partial_x^{k+2}(f(u)-\beta
		g(u_t)(\tau))\|_{L^2}\\
		&\leq&\nonumber CR^2
		\int_{t/2}^t(1+\tau)^{-\frac{n}{4}-\frac{k+1}{2}-\frac{n}{2}}d\tau\\
		&\leq&\nonumber
		CR^2(1+t)^{-\frac{n}{4}-\frac{k+1}{2}-\frac{n-1}{2}}.
	\end{eqnarray}
	Consequently we have  that
	$$\acute{J}\leq CR^2(1+t)^{-\frac{n}{4}-\frac{k+1}{2}}\eta(t).$$
	The above inequality implies
	\begin{equation}\label{equation 416}
		(1+t)^{\frac{n}{4}+\frac{k+1}{2}}\|\partial_x^k\Phi(u)_t\|_{L^2}\leq
		CE_0+CR^2.
	\end{equation}
	Combining \eqref{equation 416} and \eqref{equation 413}  and
	taking $E_0$ and $R$ suitably small, we obtain $\|\Phi(u)\|_X\leq
	R.$
	For $u,\tilde{u} \in X_R$, \eqref{equation 46} gives
	\begin{eqnarray} \nonumber
		\|\partial_x^k(\Phi(u)-\Phi(\tilde{u}))\|_{L^2}&=&
		\nonumber \int_0^t\|\partial_x^kG(t-\tau)\ast
		\Delta(f(u)-f(\tilde{u})\\
		\nonumber&&-\beta (g(u_t)-g(\tilde{u}_t))(\tau)\|_{L^2}d\tau\\
		&=&\int_0^{t/2}\|\partial_x^kG(t-\tau)\ast
		\nonumber \Delta(f(u)-f(\tilde{u})\\
		\nonumber&&-\beta (g(u_t)-g(\tilde{u}_t))(\tau)\|_{L^2}d\tau\\
		&&+\int_{t/2}^t\|\partial_x^kG(t-\tau)\ast
		\nonumber \Delta(f(u)-f(\tilde{u})\\
		\nonumber &&-\beta (g(u_t)-g(\tilde{u}_t))(\tau)\|_{L^2}d\tau\\
		\nonumber &=:&J_1+J_2
	\end{eqnarray}
	
	For the term $J_1$, we apply \eqref{equation 313} with $p=1\ ,j=0$
	and $l=0$, we arrive at
	\begin{eqnarray}
		\nonumber
		J_1&\leq& C
		\int_0^{t/2}(1+t-\tau)^{-\frac{n}{4}-\frac{k+1}{2}}\|(f(u)-f(\tilde{u})(\tau)-\beta
		(g(u_t)-g(\tilde{u}))(\tau)\|_{L^1} d\tau\\
		\nonumber&&+ C
		\int_0^{t/2}e^{-c(t-\tau)}\|\partial_x^k(f(u)-f(\tilde{u})-\beta (g(u_t-g(u_t)))(\tau)\|_{L^2}d\tau\\
		\label{equation 49}&=:&J_{11}+ J_{12}
	\end{eqnarray}
	By \eqref{equation 42}, we can estimate $J_{11}$ as
	\begin{eqnarray}
		\nonumber
		J_{11}&\leq& CR\|u-\tilde{u}\|_X
		\int_0^{t/2}(1+t-\tau)^{-\frac{n}{4}-\frac{k+1}{2}}(1+\tau)^{-\frac{n}{2}}d\tau\\
		\nonumber&\leq&
		CR\|u-\tilde{u}\|_X(1+t)^{-\frac{n}{4}-\frac{k}{2}}(1+t)^
		{-\frac{n}{4}-\frac{k}{2}}\eta(t),\nonumber
	\end{eqnarray}
	where $\eta$ be defined  in \eqref{equation 410}.
	It follows from the Gagliardo-Nirenberg inequality and
	\eqref{equation 42} that
	\begin{eqnarray}
		\nonumber
		J_{12}&\leq & \int_0^{t/2}e^{-c(t-\tau)}\Big[(\|\partial_x^k u\|_{L^2}+\|\partial_x^k\tilde{u}\|_{L^2} )\|u-\tilde{u}\|_{L^{\infty}} \\
		\nonumber &&+(\|u\|_{L^{\infty}}+\|\tilde{u}\|_
		{L^{\infty}}
		)\|\partial_x^k(u-\tilde{u})\|_{L^2}+(\|\partial_x^k u_t\|_{L^2}+\|\partial_x^k\tilde{u}_t\|_{L^2} )\|u_t-\tilde{u_t}\|_{L^{\infty}} \\
		\nonumber &&+(\|u_t\|_{L^{\infty}}+\|\tilde{u}_t\|_
		{L^{\infty}}
		)\|\partial_x^k(u_t-\tilde{u}_t)\|_{L^2}\Big] d\tau\\
		\nonumber &\leq &
		CR\int_0^{t/2}e^{-c(t-\tau)}(1+\tau)^{-\frac{n}{4}-\frac{k}{2}-\frac{n}{2}}
		\|u-\tilde{u}\|_X d\tau\\
		\nonumber &\leq &CR\|u-\tilde{u}\|_Xe^{-ct}.
	\end{eqnarray}
	Finally, we estimate term $J_2$ on the time $[t/2,t].$ Applying
	\eqref{equation 313} with $p=2$, $j=k$, $l=0$ and using
	\eqref{equation 42}, we obtain
	\begin{eqnarray}
		\nonumber J_2&\leq& C\int_{t/2}^t
		(1+t-\tau)^{-\frac{1}{2}}\|\partial_x^k(f(u)-f(\tilde{u})-\beta
		(g(u_t)-g(\tilde{u}_t))(\tau))\|_{L^2}  d\tau\\
		&&\nonumber +C\int_{t/2}^t
		e^{-c(t-\tau)}\|\partial_x^k(f(u)-f(\tilde{u})-\beta
		(g(u_t)-g(\tilde{u}_t))(\tau))\|_{L^2}\\
		\nonumber&\leq&
		CR\|u-\tilde{u}\|_X\int_{t/2}^t(1+t-\tau)^{-\frac{1}{2}}(1+t)^{-\frac{n}{4}-\frac{k}{2}-\frac{n}{2}}d\tau
		\\ \nonumber&\leq& CR\|u-\tilde{u}\|_X(1+t)^{-\frac{n}{4}-\frac{k}{2}-\frac{n-1}{2}}\label{equation
			412};
	\end{eqnarray}
	which implies
	\begin{equation}
		(1+t)^{
			\frac{n}{4}+\frac{k}{2}}\|\partial_x^k\Phi(u)-\Phi(\tilde{u})\|_{L^2}\leq CR\|u-\tilde{u}\|_X.
	\end{equation}
	Similarly, for $0\leq k\leq s$ and $u,\tilde{u}\in X$ from
	\eqref{equation 42}, \eqref{equation 314}, we deduce that
	\[\begin{split}
		\|\partial_x^k(\Phi(u)-\Phi(\tilde{u})_t\|_{L^2}&=
 \int_0^t\|\partial_x^kG_t(t-\tau)\ast
		\Delta(f(u)-f(\tilde{u}) 
 -\beta (g(u_t)-g(\tilde{u}_t))(\tau)\|_{L^2}d\tau\\
		&=\int_0^{t/2}\|\partial_x^kG_t(t-\tau)\ast
 \Delta(f(u)-f(\tilde{u}) 
	  -\beta (g(u_t)-g(\tilde{u}_t))(\tau)\|_{L^2}d\tau\\
		&\qquad+\int_{t/2}^t\|\partial_x^kG_t(t-\tau)\ast
 \Delta(f(u)-f(\tilde{u}) 
 -\beta (g(u_t)-g(\tilde{u}_t))(\tau)\|_{L^2}d\tau\\
 &=: \acute{J_1}+\acute{J_2}
	\end{split}\]
	For the term $J_1$, we use \eqref{equation 314} with $p=1,j=0$ and
	$l=2$. We have
	\[\begin{split}
		\acute{J_1}&\leq  C
		\int_0^{t/2}(1+t-\tau)^{-\frac{n}{4}-\frac{k+2}{2}}\|f(u)-f(\tilde{u})(\tau)-\beta
		(g(u_t)-g(\tilde{u}_t))(\tau)\|_{L^1} d\tau\\
	&\qquad+ C
		\int_0^{t/2}e^{-c(t-\tau)}\|\partial_x^{k+2}(f(u)-f(\tilde{u}_t)-\beta (g(u_t)-g(\tilde{u}_t))(\tau)\|_{L^2}d\tau\\
 &=:\acute{J_{11}}+ \acute{J_{12}}.
	\end{split}\]
	By \eqref{equation 42}, we have
	\[\begin{split}
		\acute{J_{11}}&\leq
		\int_0^{t/2}(1+t-\tau)^{-\frac{n}{4}-\frac{k+2}{2}}(\|u\|_{L^2}+\|\tilde{u}_{L^2})\|u-\tilde{u}\|_{L^2}\\
		&\qquad+(\|u_t\|_{L^2}+\|\tilde{u}_t\|_{L^2})(\|u_t-\tilde{u}_t\|_{L^2}) d\tau\\
 &\leq CR\|u-\tilde{u}\|_X(1+t)^{-\frac{n}{4}-\frac{k+1}{2}}
		\int_0^{t/2}(1+t-\tau)^{-\frac{n}{2}}d\tau\\
 &\leq
		CR|u-\tilde{u}\|_X(1+t)^{-\frac{n}{4}-\frac{k+1}{2}}\eta(t),
	\end{split}\]
	where $\eta$ be defined in \eqref{equation 410}.
	Also, the term $\acute{J_{12}}$ is estimated similarly  as before
	and we can estimate the term $\acute{J_{12}}$ as
	$$\acute{J_{12}}\leq CR\|u-\tilde{u}\|_X e^{-ct}.$$
	By applying  \eqref{equation 314} with $p=2$, $j=k+2$, $l=0$ and
	\eqref{equation 41}, we obtain
\[	\begin{split}
	 \acute{ J_2}&\leq 
		C\int_{t/2}^t \|\partial_x^{k+2}(f(u)-f(\tilde{u})-\beta
		(g(u_t)-g(\tilde{u}_t)(\tau))\|_{L^2}  d\tau\\
		& \qquad  +C\int_{t/2}^t
		e^{-c(t-\tau)}\|\partial_x^{k+2}(f(u)-f(\tilde{\tilde{u}})-\beta
		(g(u_t)-g(\tilde{u}_t)(\tau))\|_{L^2}d \tau\\
		&\leq  
		\int_{t/2}^t(\|\partial_x^{k+2}u\|_{L^2}+\|
		\partial_x^{k+2}\tilde{u}\|_{L^2})\|u-\tilde{u}\|_{L^{\infty}}+(\|u\|_{L^{\infty}}+\|\tilde{u}\|_{L^{\infty}})
		\|\partial_x^{k+2}(u-\tilde{u})\|_{L^2}\\
		&\qquad + ((\|\partial_x^{k+2}u_t\|_{L^2}+\|
		\partial_x^{k+2}\tilde{u}_t\|_{L^2})\|u-\tilde{u_t}\|_{L^{\infty}}+(\|u_t\|_{L^{\infty}}+\|\tilde{u_t}\|_{L^{\infty}})
		\|\partial_x^{k+2}(u-\tilde{u_t})\|_{L^2})d\tau\\
		&\leq    
		CR\|u-\tilde{u}\|_X\int_{t/2}^t(1+\tau)^{-\frac{n}{4}-\frac{k+2}{2}-\frac{n}{2}}d\tau\\
		&\leq   
		CR\|u-\tilde{u}\|_X(1+t)^{-\frac{n}{4}-\frac{k+1}{2}-\frac{n-1}{2}}\\
		&\leq 
		CR\|u-\tilde{u}\|_X(1+t)^{-\frac{n}{4}-\frac{k+1}{2}}
	\end{split}\]
	Substituting all these estimates together with the previous
	estimate and taking $R$ suitably small, yields
	\begin{equation}
		\label{equation 419}\|\Phi(u)-\Phi(\tilde{u})\|_X\leq
		\frac{1}{2}\|u-\tilde{u}\|_X.
	\end{equation}
	From \eqref{equation 419}, we deduce that $\Phi$ is strictly
	contracting mapping. Then there exists a fixed point $u\in X_R$ of
	the mapping $\Phi$, which is a solution to \eqref{equation 11},
	\eqref{equation 12}.The proof of the theorem is now complete.
\fim
The proof of  the previous theorem shows that when $n\geq 2$, the
solution $u$ to the integral equation (2.10) is asymptotic to the
linear solution $u_L(t)$ given by the formula $u_L(t)= G(t)\ast
u_1+H(t)\ast u_0$ as $t\rightarrow \infty$. This result is stated
as follows.
\begin{lemma}\label{lemma-1}
	Let $n\geq2$ and  assume the same conditions of Theorem ~(\ref{thm
		41}). Then the solution $u$of the problem \eqref{equation 11} ,
	\eqref{equation 12} which is constructed in theorem ~(\ref{thm
		41}), can be approximated by the solution $u_L$ to the linearized
	problem \eqref{equation 21}, \eqref{equation 22} as $t\rightarrow
	\infty$. More precisely, we have
	$$ \|\partial_x^k(u-u_L)(t)\|_{L^2}\leq CE_0^2(1+t)^{-\frac{n}{4}-\frac{k}{2}}\eta(t),$$
	$$ \|\partial_x^k(u-u_L)_t(t)\|_{L^2}\leq CE_0^2(1+t)^{-\frac{n}{4}-\frac{k+1}{2}}\eta(t),$$
	
	for $0\leq k\leq s+2$ and $0\leq k\leq s$, respectively, where
	$u_L(t):=G(t)\ast u _1+H(t)\ast u_0)$ is the linear solution and
	$\eta(t)$ is defined in \eqref{equation 410}.
\end{lemma}
\section{ Decay estimates of solutions for $L^2$}\label{section4}
In the previous section, we have proved global existence and
asymptotic behavior of solutions to the Cauchy problem
\eqref{equation 11}, \eqref{equation 12} with $L^1$ data.

In this section, we prove a similar decay estimate
of solution with $L^2$ data for $n\geq2$.\\
Based on the decay estimates of solutions to the linear problem
\eqref{equation 21}, \eqref{equation 22}, we define the following
solution space:
\begin{equation}
	\nonumber X=\{u \in C([0,\infty);H^{s+2}(\mathbb{R}^n))\bigcap
	C^1([0,\infty);H^s(\mathbb{R}^n)):\|u\|_X<\infty\},
\end{equation}
where
\begin{equation}
	\nonumber \|u\|_X= \sup_{t \geq 0} \left \{\sum_{k\leq
		s+2}(1+t)^{\frac{k}{2}}\|\partial_x^k u(t)\|_{L^2}+ \sum_{k\leq
		s}(1+t)^{\frac{k}{2}}\|\partial_x^k u_t(t)\|_{L^2}\right\}.
\end{equation}
For $ R>0 $, we define
$$X_R= \{u\in X:\|u\|_X\leq R\}.$$
Note that from the Gagliardo-Nirenberg inequality for $u\in X_R$,
we have
\begin{equation}
	\label{equation 51} \|u(t)\|_{L^{\infty}}\leq
	C(1+t)^{-\frac{n}{4}}.\end{equation}
\begin{Theorem}
	Suppose that $u_0 \in H^{s+2} ,u_1 \in H^s(\mathbb{R}^n)\bigcap
	\dot{W}^{-1,2}(\mathbb{R}^n)$, such that $n\geq 1,s\geq
	max\{0,[\frac{n}{2}]-1\}$, and $f(v),g(v)$ are smooth and
	satisfies $f(v)= O(v^2),g(v)=O(v^2)$ for $v\rightarrow 0$. Let
	\begin{equation}
		\nonumber E_1:= \|u_0\|_{L^2}+
		\|u_1\|_{\dot{W}^{-1,2}}+\|u_0\|_{H^{s+2}}+ \|u_1\|_{H^s}.
	\end{equation}
		If $E_0$ is suitably small, the Cauchy problem \eqref{equation 11}
	and \eqref{equation 12} has a unique global solution $u(x,t)$
	satisfying
	\begin{equation}
		\nonumber X=u \in C([0,\infty);H^{s+2}(\mathbb{R}^n))\bigcap
		C^1([0,\infty);H^s(\mathbb{R}^n)).
	\end{equation}
The solution $u$ also satisfies the decay estimate
	\begin{equation}
		\label{equation 52}
		\|\partial_x^k u(t)\|_{L^2}\leq C E_1(1+t)^{-\frac{k}{2}}
	\end{equation}
	and
	\begin{equation}
		\label{equation 53}
		\|\partial_x^h u_t(t)\|_{L^2}\leq C E_1(1+t)^{-\frac{h+1}{2}}.
	\end{equation}
	for $0\leq k\leq s+2$ and $0\leq h \leq s.$
	\end{Theorem}
\proof
	Let the  mapping $\Phi$ be defined in \eqref{equation 46}.
	Applying $\partial_x^k$ to $\Phi$ and take $L^2$ norm. We have
	\begin{eqnarray}\label{equation 54}
		\nonumber
		\|\partial_x^k\Phi(u)\|_{L^2}&\leq &\|\partial_x^k G(t)\ast
		u_1\|_{L^2}+ \|\partial_x^kH(t)\ast u_0\|_{L^2}\\
		&&
		\nonumber+C\int_0^t\|\partial_x^kG(t-\tau)\ast
		\Delta(f(u)-\beta g(u_t))(\tau)\|_{L^2}d\tau\\
		&& :=I_1+I_2+J.
	\end{eqnarray}
	We use \eqref{equation 39} with $p=2,j=l=0$ and get
	\begin{eqnarray}
		\label{equation 55}
		I_1&\leq & C (1+t)^{-\frac{k}{2}}
		\|u_1\|_{\dot{W}^{-1,2}}+ C
		e^{-ct}\|\partial_x^{(k-2)_+}u_1\|_{L^2}\\
		\nonumber& \leq &C E_1 (1+t)^{-  \frac{k}{2}},
	\end{eqnarray}
	where $(k-2)_+= \max\{k-2,0\}.$
	By  applying \eqref{equation 310} with $p=2,j=l=0$, we get
	\begin{eqnarray}
		\label{equation{56}}
		I_2&\leq &C (1+t)^{-\frac{k}{2}}
		\|u_0\|_{L^2}+ C e^{-ct}\|\partial_x^{k}u_0\|_{H^{s+2}}\\
		\nonumber & \leq & C E_1 (1+t)^{- \frac{k}{2}}.
	\end{eqnarray}
	To estimate the nonlinear $J$, as in the pervious section, we
	divide as $J=J_1+J_2$ where $J_1$ and $J_2$ correspond to the time
	intervals$[0,t/2]$ and $[t/2,t]$, respectively. For the term
	$J_1$, we use \eqref{equation 313} with $p=1,j=l=0$ and deduce
	that
	\begin{eqnarray}
		\nonumber
		J_1&\leq& C
		\int_0^{t/2}(1+t-\tau)^{-\frac{n}{4}-\frac{k+1}{2}}\|f(u)(\tau)-\beta
		g(u_t)(\tau)\|_{L^1} d\tau\\
		\nonumber&&+ C
		\int_0^{t/2}e^{-c(t-\tau)}\|\partial_x^k(f(u)-\beta g(u_t))(\tau)\|_{L^2}d\tau\\
		\label{equation 57}&=:&J_{11}+ J_{12}
	\end{eqnarray}
	By \eqref{equation 41}, we have $\|f(u)-\beta
	g(u_t)(\tau)\|_{L^1}\leq CR^2$. Thus we can estimate the $J_{11}$
	as
	\begin{eqnarray}
		\nonumber
		J_{11}&\leq& CR^2
		\int_0^{t/2}(1+t-\tau)^{-\frac{n}{4}-\frac{k+1}{2}}(1+\tau)^{-\frac{n}{2}}d\tau\\
		\nonumber&\leq&
		CR^2(1+t)^{-\frac{k}{2}}\int_0^{t/2}(1+\tau)^{-\frac{n}{4}-\frac{1}{2}}d\tau\\
		&\leq&CR^2(1+t)^{-\frac{k}{2}}\nonumber.
	\end{eqnarray}
	By applying \eqref{equation 41}  and Gagliardo-Nirenberg
	inequality, we get
	\begin{eqnarray}
		\nonumber\|\partial_x^k f(u)\|_{L^2}&\leq& C
		\|u\|_{L^{\infty}}\|\partial_x^ku\|_{L^2}\leq
		C(1+t)^{\frac{n}{4}-\frac{k}{2}}R^2\\
		|\partial_x^k g(u_t)\|_{L^2}&\leq& C
		\|u_t\|_{L^{\infty}}\|\partial_x^ku_t\|_{L^2}\leq
		C(1+t)^{\frac{n}{4}-\frac{k}{2}}R^2\label{equation 58}
	\end{eqnarray}
	Thus we have $$J_{12} \leq
	CR^2\int_0^{t/2}e^{-c(t-\tau)}(1+\tau)^{-\frac{n}{4}-\frac{k}{2}}d\tau
	\leq CR^2e^{-ct}. $$
	
	It follows from \eqref{equation 313} with $p=1,j=k$ and $l=2$ that
	\begin{eqnarray}
		\nonumber J_2&\leq& C\int_{t/2}^t
		(1+t-\tau)^{-\frac{n}{4}-\frac{1}{2}}\|\partial_x^k(f(u)-\beta
		g(u_t))(\tau)\|_{L^1}  d\tau\\
		&&\nonumber +C\int_{t/2}^t
		e^{-c(t-\tau)}\|\partial_x^{k+2}(f(u)-\beta
		g(u_t))(\tau)\|_{L^2}\\
		\nonumber &&=:J_{21}+J_{22}
	\end{eqnarray}
We have using \eqref{equation 42} that
	\begin{eqnarray}
		\nonumber J_{21}&\leq& C\int_{t/2}^t
		(1+t-\tau)^{-\frac{n}{4}-\frac{1}{2}}(\|u\|_{L^2}\|
		\|\partial_x^ku\|_{L^2}+\|u_t\|_{L^2}\|\partial_x^ku_t\|_{L^2} )
		d\tau\\
		\nonumber &\leq& CR^2\int_{t/2}^t
		(1+t-\tau)^{-\frac{n}{4}-\frac{1}{2}}(1+\tau)^{-\frac{k}{2}}d\tau\\
		\nonumber&\leq& C R^2(1+t)^{-\frac{k}{2}}\int_{t/2}^t
		(1+t-\tau)^{-\frac{n}{4}-\frac{1}{2}}d\tau\\
		\nonumber&\leq&CR^2(1+t)^{-\frac{k}{2}}
	\end{eqnarray}
	To estimate the term $J_{22}$, by \eqref{equation 58}, we have
	\begin{eqnarray}
		\nonumber J_{22}&\leq& C\int_{t/2}^t
		e^{-c(t-\tau)}\|\partial_x^{k+2}(f(u)-\beta
		g(u_t))(\tau)\|_{L^2}d\tau\\
		&\leq&
		\nonumber CR^2\int_{t/2}^te^{-c({t-\tau})}(1+\tau)^{-\frac{n}{4}-\frac{k+2}{2}} d\tau\\
		\nonumber&\leq&CR^2(1+t)^{-\frac{k}{2}}.
	\end{eqnarray}
	
	The above  inequality shows that
	\begin{equation}
		\label{equation 510} (1+t)^{\frac{k}{2}}
		\|\partial_x^k\Phi(u)\|\leq CE_1+CR^2
	\end{equation}
	We deduce from \eqref{equation 46} that
	\begin{equation}
		\label{equation 511} \Phi(u)_t= G_t(t)\ast u_1+ H_t(t)\ast
		u_0+\int_0^tG_t(t-\tau)\ast \Delta(f(u)-\beta
		g(u_t))(\tau)\|_{L^2}d\tau.
	\end{equation}
	Applying $\partial _x^k$ to $\Phi(u)_t$ and taking $L^2$-norm we have
	\begin{eqnarray}\label{equation 512}
		\nonumber
		\|\partial_x^k\Phi(u)_t\|_{L^2}&\leq &\|\partial_x^k G_t(t)\ast
		u_1\|_{L^2}+ \|\partial_x^kH_t(t)\ast u_0\|_{L^2}\\
		&&
		\nonumber+C\int_0^t\|\partial_x^kG_t(t-\tau)\ast
		\Delta(f(u)-\beta g(u_t))(\tau)\|_{L^2}d\tau\\
		&& =:\acute{I_1}+ \acute{I_2}+\acute{J},
	\end{eqnarray}
	
	To estimate the term $\acute{I_1}$, apply \eqref{equation 311} with
	$p=2,l=j=0$. It yields
	\begin{eqnarray} \nonumber
		\acute{I_1}\leq C (1+t)^{-\frac{k+1}{2}}
		\|u_1\|_{\dot{W}^{-1,2}}+ C e^{-ct}\|\partial_x^{k}u_1\|_{L^2}
		\leq C E_1 (1+t)^{- \frac{k+1}{2}}.
	\end{eqnarray}
	Similarly, using \eqref{equation 312} with $p=2,j=l=0$, we have
	\begin{eqnarray}
		\nonumber
		\acute{I_2}\leq C (1+t)^{-\frac{k+1}{2}}
		\|u_0\|_{L^2}+ C e^{-ct}\|\partial_x^{k+2}u_0\|_{L^2} \leq C E_1
		(1+t)^{- \frac{k+1}{2}}.
	\end{eqnarray}
	To estimate the nonlinear term $\acute{J}$,let
\[	\begin{split}
	 \acute{J}
		&=
	  C\int_0^{t/2}\|\partial_x^kG_t(t-\tau)\ast
		\Delta(f(u)-\beta g(u_t))(\tau)\|_{L^2}d\tau\\
		&\qquad +C\int_{t/2}^{t}\|\partial_x^kG_t(t-\tau)\ast
		\Delta(f(u)-\beta g(u_t))(\tau)\|_{L^2}d\tau\\
  & =: \acute{I_1}+ \acute{I_2}+\acute{J},
	\end{split}\]
	Using \eqref{equation 314} with $p=1,j=l=0$, it yields
	\begin{eqnarray}
		\nonumber
		\acute{J_1}&\leq& C
		\int_0^{t/2}(1+t-\tau)^{-\frac{n}{4}-\frac{k+2}{2}}\|(f(u)(\tau)-\beta
		g(u_t))(\tau)\|_{L^1} d\tau\\
		\nonumber&+& C
		\int_0^{t/2}e^{-c(t-\tau)}\|\partial_x^{k+2}(f(u)-\beta g(u_t))(\tau)\|_{L^2}d\tau\\
		\nonumber&:=&\acute{J_{11}}+ \acute{J_{12}}.
	\end{eqnarray}
	By \eqref{equation 41}, we  obtain
	$$\|(f(u)-\beta g(u_t))(\tau)\|_{L^1}\leq
	CR^2.$$
	Thus we can estimate $\acute{J_{11}}$ as
\[	\begin{split}
		\acute{J_{11}}&\leq CR^2
		\int_0^{t/2}(1+t-\tau)^{-\frac{n}{4}-\frac{k+2}{2}}d\tau\\
&\leq CR^2(1+t)^{-\frac{k+1}{2}}\int_0^{t/2}(1+t-\tau)^{-\frac{n}{4}-\frac{1}{2}}d\tau\\
&\leq CR^2(1+t)^{-\frac{k+1}{2}}
	\end{split}\]
	For the term $\acute{J_{12}}$,  we have from \eqref{equation 58} that
	\begin{eqnarray}
		\nonumber
		\acute{J_{12}}&\leq& CR^2
		\int_0^{t/2}e^{-c(t-\tau)}(1+t)^{-\frac{n}{4}-\frac{k+2}{2}}d\tau\\
		\nonumber&\leq&CR^2e^{-ct}
	\end{eqnarray}
	Applying \eqref{equation 314} with $p=2,j=k+2$ and $l=0$, we get
	\begin{eqnarray}
		\nonumber\acute{ J_2}&\leq& C\int_{t/2}^t
		\|\partial_x^{k+2}(f(u)-\beta
		g(u_t)(\tau))\|_{L^2}  d\tau\\
		&&\nonumber +C\int_{t/2}^t
		e^{-c(t-\tau)}\|\partial_x^{k+2}(f(u)-\beta
		g(u_t)(\tau))\|_{L^2}\\
		&\leq&\nonumber CR^2
		\int_{t/2}^t(1+\tau)^{-\frac{n}{4}-\frac{k+2}{2}}d\tau\\
		&\leq&\nonumber
		CR^2(1+t)^{-\frac{k+1}{2}}\int_{t/2}^t(1+\tau)^{-\frac{n}{4}-\frac{1}{2}}d\tau\\
		&\leq&\nonumber
		CR^2(1+t)^{-\frac{k+1}{2}}.
	\end{eqnarray}
	Thus we have
	
	\begin{equation}\label{equation 513}
		(1+t)^{-\frac{k+1}{2}}\|\partial_x^k\Phi(u)_t\|_{L^2}\leq
		CE_1+CR^2.
	\end{equation}
	Combining \eqref{equation 510} and \eqref{equation 513} and taking
	$E_0$ and$R$ suitably small, we obtain $\|\Phi(u)\|_X\leq R.$\\
	
	For $u,\tilde{u} \in X_R$, by using \eqref{equation 46} we obtain
	\begin{eqnarray} \nonumber
		\|\partial_x^k(\Phi(u)-\Phi(\tilde{u}))\|_{L^2}&=&
		\nonumber \int_0^t\|\partial_x^kG(t-\tau)\ast
		\Delta(f(u)-f(\tilde{u})-\beta (g(u_t))-g(\tilde{u}_t))(\tau)\|_{L^2}d\tau\\
		&=&\int_0^{t/2}\|\partial_x^kG(t-\tau)\ast
		\nonumber \Delta(f(u)-f(\tilde{u})-\beta (g(u_t))-g(\tilde{u}_t))(\tau)\|_{L^2}d\tau\\
		&&+\int_{t/2}^t\|\partial_x^kG(t-\tau)\ast
		\nonumber \Delta(f(u)-f(\tilde{u})-\beta (g(u_t))-g(\tilde{u}_t))(\tau)\|_{L^2}d\tau\\
		\nonumber &=:&J_1+J_2
	\end{eqnarray}
	By applying \eqref{equation 313} with $p=1,j=0$ and $l=0$, we have
\[	\begin{split}
			J_1&\leq C
		\int_0^{t/2}(1+t-\tau)^{-\frac{n}{4}-\frac{k+1}{2}}\|(f(u)-f(\tilde{u})(\tau)-\beta
		(g(u_t)-g(\tilde{u_t}))(\tau)\|_{L^1} d\tau\\
&\qquad+ C
		\int_0^{t/2}e^{-c(t-\tau)}\|\partial_x^k(f(u)-f(\tilde{u})-\beta (g(u_t-g(\tilde{u}_t)))(\tau)\|_{L^2}d\tau\\
 &=:J_{11}+ J_{12}
	\end{split}\]
	Using \eqref{equation 42}, we get
	\begin{eqnarray}
		\nonumber
		J_{11}&\leq& CR\|u-\tilde{u}\|_X
		\int_0^{t/2}(1+t-\tau)^{-\frac{n}{4}-\frac{k+1}{2}}d\tau\\
		\nonumber&\leq&
		CR\|u-\tilde{u}\|_X(1+t)^{-\frac{k}{2}}\int_0^{t/2}(1+t-\tau)^
		{-\frac{n}{4}-\frac{1}{2}}d\tau.\nonumber\\
		\nonumber&\leq&CR\|u-\tilde{u}\|_X(1+t)^{-\frac{k}{2}}
	\end{eqnarray}
	Also, we have
	\begin{eqnarray} \nonumber
		J_{12}&\leq & \int_0^{t/2}e^{-c(t-\tau)}\Big[(\|\partial_x^k u\|_{L^2}+\|\partial_x^k\tilde{u}\|_{L^2} )\|u-\tilde{u}\|_{L^{\infty}} \\
		\nonumber &&+(\|u\|_{L^{\infty}}+\|\tilde{u}\|_
		{L^{\infty}}
		)\|\partial_x^k(u-\tilde{u})\|_{L^2}+(\|\partial_x^k u_t\|_{L^2}+\|\partial_x^k\tilde{u}_t\|_{L^2} )\|u_t-\tilde{u_t}\|_{L^{\infty}} \\
		\nonumber &&+(\|u_t\|_{L^{\infty}}+\|\tilde{u}_t\|_
		{L^{\infty}}
		)\|\partial_x^k(u_t-\tilde{u}_t)\|_{L^2}\Big] d\tau\\
		\nonumber &\leq &
		CR\|u-\tilde{u}\|_X\int_0^{t/2}e^{-c(t-\tau)}(1+\tau)^{-\frac{n}{4}-\frac{k}{2}}
		d\tau\\
		\nonumber &\leq &CR\|u-\tilde{u}\|_Xe^{-ct}.
	\end{eqnarray}
	To estimate the term $J_2$, apply \eqref{equation 313} with
	$p=1,j=k,l=2.$ We obtain
	\begin{eqnarray}
		\nonumber J_2&\leq& C\int_{t/2}^t
		(1+t-\tau)^{-\frac{n}{4}-\frac{1}{2}}\|\partial_x^k(f(u)-f(\tilde{u})-\beta
		(g(u_t)-g(\tilde{u}_t))(\tau))\|_{L^1}  d\tau\\
		&&\nonumber +C\int_{t/2}^t
		e^{-c(t-\tau)}\|\partial_x^{k+2}(f(u)-f(\tilde{u})-\beta
		(g(u_t)-g(\tilde{u}_t))(\tau))\|_{L^2}d\tau\\
		\nonumber&=:&J_{21}+J_{22}
	\end{eqnarray}
	By using \eqref{equation 42}, we get
	\begin{eqnarray}
		\nonumber
		J_{21}&\leq & \int_{t/2}^t(1+t-\tau)^{-\frac{n}{4}-\frac{1}{2}}\Big[(\|\partial_x^k u\|_{L^2}+\|\partial_x^k\tilde{u}\|_{L^2} )\|u-\tilde{u}\|_{L^2} \\
		\nonumber &&+(\|u\|_{L^{2}}+\|\tilde{u}\|_
		{L^2})
		)\|\partial_x^k(u-\tilde{u})\|_{L^2}+(\|\partial_x^k u_t\|_{L^2}+\|\partial_x^k\tilde{u}_t\|_{L^2} )\|u_t-\tilde{u_t}\|_{L^2} \\
		\nonumber &&+(\|u_t\|_{L^2}+\|\tilde{u}_t\|_
		{L^2}
		)\|\partial_x^k(u_t-\tilde{u}_t)\|_{L^2}\Big] d\tau\\
		\nonumber &\leq&
		CR\|u-\tilde{u}\|_X\int_{t/2}^{t}(1+t-\tau)^{-\frac{n}{4}-\frac{1}{2}}(1+\tau)^{-\frac{k}{2}}
		d\tau\\
		\nonumber &\leq&CR\|u-\tilde{u}\|_X(1+t)^{-\frac{k}{2}}
	\end{eqnarray}
	
	Finally, we estimate the term $J_{22}$ as
	\begin{eqnarray} \nonumber
		J_{22}&\leq & \int_{t/2}^{t}e^{-c(t-\tau)}\Big[(\|\partial_x^{k+2} u\|_{L^2}+\|\partial_x^{k+2}\tilde{u}\|_{L^2} )\|u-\tilde{u}\|_{L^2} \\
		\nonumber &&+(\|u\|_{L^2}+\|\tilde{u}\|_
		{L^2}
		)\|\partial_x^k(u-\tilde{u})\|_{L^2}+(\|\partial_x^k u_t\|_{L^2}+\|\partial_x^k\tilde{u}_t\|_{L^2} )\|u_t-\tilde{u_t}\|_{L^2} \\
		\nonumber &&+(\|u_t\|_{L^2}+\|\tilde{u}_t\|_
		{L^2}
		)\|\partial_x^k(u_t-\tilde{u}_t)\|_{L^2}\Big] d\tau\\
		\nonumber &\leq &
		CR\|u-\tilde{u}\|_X\int_{t/2}^{t}e^{-c(t-\tau)}(1+\tau)^{-\frac{n}{4}-\frac{k+2}{2}}
		d\tau\\
		\nonumber &\leq &CR\|u-\tilde{u}\|_X(1+\tau)^{\frac{k}{2}}
	\end{eqnarray}
	Thus we have shown that
	\begin{equation}
		\label{quation513}(1+t)^{
			\frac{k}{2}}\|\partial_x^k(\Phi(u)-\Phi(\tilde{u}))\|_{L^2}\leq
		CR\|u-\tilde{u}\|_X.
	\end{equation}
	Suppose that $u,\tilde{u} \in X_{R}.$ It follows from
	\eqref{equation 46} that
	\begin{eqnarray} \nonumber
		\|\partial_x^k(\Phi(u)-\Phi(\tilde{u})_t\|_{L^2}&=&
		\nonumber \int_0^t\|\partial_x^kG_t(t-\tau)\ast
		\Delta(f(u)-f(\tilde{u})-\beta (g(u_t))-g(\tilde{u}_t))(\tau)\|_{L^2}d\tau\\
		&=&\int_0^{t/2}\|\partial_x^kG_t(t-\tau)\ast
		\nonumber \Delta(f(u)-f(\tilde{u})-\beta (g(u_t))-g(\tilde{u}_t))(\tau)\|_{L^2}d\tau\\
		&&+\int_{t/2}^t\|\partial_x^kG_t(t-\tau)\ast
		\nonumber \Delta(f(u)-f(\tilde{u})-\beta (g(u_t))-g(\tilde{u}_t))(\tau)\|_{L^2}d\tau\\
		\nonumber &=:&\acute{J_1}+\acute{J_2}.
	\end{eqnarray}
	By using \eqref{equation 314} with $p=1,j=0$, we have
	\begin{eqnarray} \nonumber
		\acute{J_1}&\leq& C
		\int_0^{t/2}(1+t-\tau)^{-\frac{n}{4}-\frac{k+2}{2}}\|f(u)-f(\tilde{u})(\tau)-\beta
		(g(u_t)-g(\tilde{u}_t))(\tau)\|_{L^1} d\tau\\
		\nonumber&+& C
		\int_0^{t/2}e^{-c(t-\tau)}\|\partial_x^{k+2}(f(u)-f(\tilde{u}_t)-\beta (g(u_t)-g(\tilde{u}_t))(\tau)\|_{L^2}d\tau\\
		\nonumber&=:&\acute{J_{11}}+ \acute{J_{12}}.
	\end{eqnarray}
	By \eqref{equation 42}, we obtain
	\begin{eqnarray}
		\nonumber
		\acute{J_{11}}&\leq&
		\int_0^{t/2}(1+t-\tau)^{-\frac{n}{4}-\frac{k+2}{2}}(\|u\|_{L^2}+\|\tilde{u}_{L^2})\|u-\tilde{u}\|_{L^2}\\
		&&\nonumber\times(\|u_t\|_{L^2}+\|\tilde{u}_t\|_{L^2})(\|u_t-\tilde{u}_t\|_{L^2}) d\tau\\
		\nonumber &\leq &CR\|u-\tilde{u}\|_X(1+t)^{-\frac{k+1}{2}}
		\int_0^{t/2}(1+t\tau)^{-\frac{n}{4}-\frac{1}{2}}d\tau\\
		\nonumber &\leq &CR|u-\tilde{u}\|_X(1+t)^{-\frac{k+1}{2}}.
	\end{eqnarray}
	For the term $\acute{J_{12}}$, by \eqref{equation 42} we get
\[	\begin{split}
\acute{J_{12}}&\leq
\int_0^{t/2}e^{-c(t-\tau)}\left[(\|\partial_x^{k+2} u\|_{L^2}+\|\partial_x^{k+2}\tilde{u}\|_{L^2} )\|u-\tilde{u}\|_{L^{\infty}}\right. \\
&\qquad+(\|u\|_{L^{\infty}}+\|\tilde{u}\|_
{L^{\infty}})
)\|\partial_x^{k+2}(u-\tilde{u})\|_{L^2}+(\|\partial_x^{k+2} u_t\|_{L^2}+\|\partial_x^{k+2}\tilde{u}_t\|_{L^2} )\|u_t-\tilde{u_t}\|_{L^{\infty}} \\
&\qquad\left.+(\|u_t\|_{L^{\infty}}+\|\tilde{u}_t\|_
{L^{\infty}}
)\|\partial_x^{k+2}(u_t-\tilde{u}_t)\|_{L^2}\right] d\tau\\
&\leq CR\|u-\tilde{u}\|_X\int_0^{t/2}e^{-c({t-\tau})}(1+\tau)^{-\frac{k+2}{2}}
(1+\tau)^{-\frac{n}{4}}d\tau\\
&\leq  CR\|u-\tilde{u}\|_Xe^{-ct}
\end{split}\]
	By applying \eqref{equation 314} with $p=2$, $j=k+2$ and $l=0$, we
	conclude that
\[	\begin{split}
	\acute{ J_2}&\leq
		C\int_{t/2}^t \|\partial_x^{k+2}(f(u)-f(\tilde{u})-\beta
		(g(u_t)-g(\tilde{u}_t)(\tau))\|_{L^2}  d\tau\\
		&\qquad  +C\int_{t/2}^t
		e^{-c(t-\tau)}\|\partial_x^{k+2}(f(u)-f(\tilde{\tilde{u}})-\beta
		(g(u_t)-g(\tilde{u}_t)(\tau))\|_{L^2}d \tau\\
		&\leq 
		\int_{t/2}^t(\|\partial_x^{k+2}u\|_{L^2}+\|
		\partial_x^{k+2}\tilde{u}\|_{L^2})\|u-\tilde{u}\|_{L^{\infty}}+(\|u\|_{L^{\infty}}+\|\tilde{u}\|_{L^{\infty}})
		\|\partial_x^{k+2}(u-\tilde{u})\|_{L^2}\\
		&\qquad+ (\|\partial_x^{k+2}u_t\|_{L^2}+\|
		\partial_x^{k+2}\tilde{u}_t\|_{L^2})\|u-\tilde{u_t}\|_{L^{\infty}}+(\|u_t\|_{L^{\infty}}+\|\tilde{u_t}\|_{L^{\infty}})
		\|\partial_x^{k+2}(u-\tilde{u_t})\|_{L^2})d\tau\\
		&\leq  
		CR\|u-\tilde{u}\|_X\int_{t/2}^t(1+\tau)^{-\frac{n}{4}-\frac{k+2}{2}}d\tau\\
		&\leq  
		CR\|u-\tilde{u}\|_X(1+t)^{-\frac{k+1}{2}}\int_{t/2}^t(1+\tau)^{-\frac{n}{4}-\frac{1}{2}}d\tau\\
		&\leq  
		CR\|u-\tilde{u}\|_X(1+t)^{-\frac{k+1}{2}}.
	\end{split}\]
	Consequently, we have shown that
	\begin{equation}	\label{equation514}
	(1+t)^{\frac{k+1}{2}}\|\
		\partial_x^k( \Phi(u)-\Phi(\tilde{u}))_t\|_X\leq
		CR\|u-\tilde{u}\|_X.
	\end{equation}
	Using \eqref{equation 513} and \eqref{equation514} and taking $R$
	suitably small, it yields
	\begin{equation}
		\label{equation 515}\|\Phi(u)-\Phi(\tilde{u})\|_X\leq
		\frac{1}{2}\|u-\tilde{u}\|_X.
	\end{equation}
	From \eqref{equation 515}, we conclude that $\Phi$ is a
	contracting mapping. Then there exists a fixed point $u\in X_R$ of
	mapping $\Phi$, which is a solution  \eqref{equation 11} and \eqref{equation 12} and the proof is
	completed.
	\fim
Finally we   study the asymptotic linear profile of the solution.

Suppose that $u_L$ given by the formula $u_L(t)=G(t)\ast
u_1+H(t)\ast u_0$. In the previous two section, we have shown that
the solution $u$ to the problem \eqref{equation 11} and \eqref{equation 12} can be approximated by
the linear solution $u_L.$ Now the aim  is to derive
a simpler asymptotic profile of the linear solution $u_L.$

In the Fourier space, we obtain  $\hat{u}_L(\xi,t)=
\hat{G}(\xi,t)\hat{u}_1(\xi,t)+\hat{H}(\xi,t)\hat{u}_0(\xi)$,
where $\hat{G}(\xi,t)$ and $\hat{H}(\xi,t)$ are given explicitly
in \eqref{equation 27} and \eqref{equation 28}. First we give the asymptotic expansions of
$\hat{G}(\xi,t)$ and $\hat{H}(\xi,t)$ for $\xi\rightarrow 0$. By
Using the Taylor expansion to \eqref{equation 25}, we obtain
\begin{eqnarray}
	\nonumber
	\lambda_{\pm}(\xi)&=&\frac{1}{2}(\alpha|\xi|^2-|\xi|^4)\pm
	\frac{|\xi| i}{2}(2+|\xi|^2-\frac{\alpha^2}{4}|\xi|^4+O(|\xi|^4)\\
	\nonumber &=&\pm i|\xi|+\frac{\alpha}{2}|\xi|^2+O(|\xi|^3)
\end{eqnarray}
and
\begin{eqnarray}
	\nonumber
	\frac{1}{\lambda_+-\lambda_-}&=&\frac{1}{i|\xi|\sqrt{4+4|\xi|^2-|\xi|^6-\alpha
			^2|\xi|^4+2\alpha|\xi|^5}}\\
	\nonumber&=&\frac{1}{2i|\xi|}(1-\frac{1}{2}|\xi|^2+O(|\xi|^4)).
\end{eqnarray}
Substituting these expansions to \eqref{equation 27} and \eqref{equation 28}, we obtain
\begin{eqnarray}
	\nonumber
	\hat{G}(\xi,t)&=&\frac{e^{\lambda_+t}-e^{\lambda_-t}}{\lambda_+-\lambda_-}\\
	\nonumber&=&\frac{1}{2i|\xi|}
	\Big(e^{\frac{\alpha}{2}|\xi|^2t}(e^{|\xi|ti}-e^{-|\xi|ti})+e^{\frac{\alpha}{2}|\xi|^2t}(O(|\xi|^2)+O(|\xi|^3t))\Big)
\end{eqnarray}
and
\begin{eqnarray}
	\nonumber
	\hat{H}(\xi,t)&=&\frac{\lambda_+e^{\lambda_-t}-\lambda _-e^{\lambda_+t}}{\lambda_+-\lambda_-}\\
	\nonumber&=&\frac{1}{2}
	e^{\frac{\alpha}{2}|\xi|^2t}(e^{|\xi|ti}+e^{-|\xi|ti})+e^{\frac{\alpha}{2}|\xi|^2t}(O(|\xi|^2)+O(|\xi|^3t))
\end{eqnarray}
for $\xi\rightarrow 0$. Let
$$\hat{G_0}(\xi,t)=\frac{1}{2i|\xi|}e^{\frac{\alpha}{2}|\xi|^2t}(e^{i|\xi|t}-e^{-i|\xi|t})$$
and
\[
\hat{H_0}(\xi,t)=\frac{1}{2}e^{\frac{\alpha}{2}|\xi|^2t}(e^{i|\xi|t}+e^{-i|\xi|t}).
\]
Thus for $|\xi|\leq r_0$ we obtain
$$|(\hat{G}-\hat{G_0})(\xi,t)|\leq Ce^{-c|\xi|^2t},\ \ \ \ \ |(\hat{H}-\hat{H_0})(\xi,t)|\leq C|\xi|e^{-c|\xi|^2t}$$
where $r_0$ is a small positive constant. Now we define
$\overline{u}_L$ by
\begin{equation}
	\label{equation 61} \overline{u}_L(t)=G_0(t)\ast u_1+H_0(t)\ast
	u_0.
\end{equation}
$\overline{u}_L$
gives an asymptotic profile of the linear solution $u_L$.  
\begin{Theorem}
	Suppose that $n\geq 1$ , $s\geq 0$ and $u_0 \in H^{s+2}\bigcap
	L^1$ and $u_1 \in H^s\bigcap \dot{W}^{-1,1}$. Put
	$E_0=\|u_0\|_{L^1}+\|u_1\|_{\dot{W}^{-1,1}}+\|u_0\|_{H^{s+2}}+\|u_1\|_{H^s}$. Let
	$u_L$ be the linear solution and $\overline{u}_L$ be defined by
	\eqref{equation 61}. Thus we have
	\begin{equation}
		\|\partial_x^k(u_L-\overline{u}_L)(t)\|_{L^2}\leq
		CE_0(1+t)^{-\frac{n}{4}-\frac{k+1}{2}}
	\end{equation}
	for $0\leq k\leq s+2.$
\end{Theorem}
\proof
	It follows from definition that
	$$(u_L-\overline{u}_L)(t)= (G-G_0)(t)\ast u_1+(H-H_0)(t)\ast u_0.$$
So  it suffices to show the following
	estimates:
	\begin{align}
		\nonumber
		\|\partial_x^k( G-G_0)(t)\ast u_1\|_{L^2} \leq C(1+t)^{-\frac{n}{2}(\frac{1}{p}-\frac{1}{2})-
			\frac{k+1-j}{2}}\ \|\partial_x^j u_1\|_{\dot{W}^{-1 ,p}} +C
		e^{-ct} \|\partial_x^{k+l-2}u_1\|_{L^2},
	\end{align}
	\begin{equation}
		\nonumber
		\|\partial_x^k( H-H_0)(t)\ast u_0\|_{L^2} \leq C(1+t)^{-\frac{n}{2}(\frac{1}{p}-\frac{1}{2})-
			\frac{k+1-j}{2}}\ \|\partial_x^j u_0\|_{L^p} +C
		e^{-ct} \|\partial_x^{k+l}\phi\|_{L^2},
	\end{equation}
	where $1\leq p\leq 2$, and $k,j$ and $l$ are nonnegative integers
	such that $0\leq j\leq k+1$. We assumed $k+l-2\geq 0$ in the first
	estimate. These estimates can be proved similarly as in the proof
	of Lemma \ref{lemma-1} by using \eqref{equation 52} for $|\xi|\leq r_0$ and \eqref{equation 37} and  \eqref{equation 37} and
	\eqref{equation 511} for $|\xi|\geq r_0$. We omit the details.
\fim


\section*{}

\end{document}